\documentclass[amstex,12pt,russian,amssymb]{article}

\usepackage{mathtext}
\usepackage[cp1251]{inputenc}
\usepackage[T2A]{fontenc}
\usepackage[russian]{babel}
\usepackage[dvips]{graphicx}
\usepackage{amsmath}
\usepackage{amssymb}
\usepackage{amsxtra}
\usepackage{latexsym}
\usepackage{ifthen}

\begin{document}

\newcounter{lemma}
\newcommand{\lemma}{\par \refstepcounter{lemma}%
{\bf Лемма \arabic{lemma}.}}

\newcounter{corollary}
\newcommand{\corollary}{\par \refstepcounter{corollary}%
{\bf Следствие \arabic{corollary}.}}

\newcounter{remark}
\newcommand{\remark}{\par \refstepcounter{remark}%
{\bf Замечание \arabic{remark}.}}

\newcounter{theorem}
\newcommand{\theorem}{\par \refstepcounter{theorem}%
{\bf Теорема \arabic{theorem}.}}

\newcounter{proposition}
\newcommand{\proposition}{\par \refstepcounter{proposition}%
{\bf Предложение \arabic{proposition}.}}

\renewcommand{\refname}{\centerline{\bf Список литературы}}

\newcommand{\proof}{{\it Доказательство.\,\,}}

\noindent УДК 517.5

\medskip
{\bf Е.А..~Севостьянов, С.А.~Скворцов} (Житомирский государственный
университет им.\ И.~Франко)

\medskip
{\bf E.A.~Sevost'yanov, S.A.~Skvortsov} (Zhitomir State University
of I.~Franko)

\medskip\medskip
{\bf О равностепенной непрерывности обобщённых квазиизометрий на
римановых многообразиях}

\medskip\medskip
{\bf On equicontinuity of generalized quasiisometries on Riemannian
manifolds}

\medskip\medskip
Настоящая работа посвящена изучению отображений с конечным
искажением на римановых многообразиях. Доказаны теоремы о локальном
поведении обобщённых квазиизометрий с неограниченной
характеристикой. В частности, нами доказано, что семейство
отображений $f:D\rightarrow {\Bbb M}_*^n$ между римановыми
многообразиями ${\Bbb M}^n$ и ${\Bbb M}_*^n$ является равностепенно
непрерывным, как только значения отображения $f$ лежат в шаре $B_R,$
$f$ не принимает значения из фиксированного континуума $K\subset
B_R,$ а характеристика квазиконформности $Q(x)$ имеет конечное
среднее колебание в каждой точке.

\medskip\medskip
The present paper is devoted to the study of mappings with finite
distortion on Riemannian manifolds. Theorems on local behavior of
generalized quasiisometries with unbounded cha\-rac\-te\-ris\-tic of
qua\-si\-con\-for\-ma\-li\-ty are obtained. In particular, we have
proved that a family of mappings $f:D\rightarrow {\Bbb M}_*^n$
between Riemannian manifolds ${\Bbb M}^n$ and ${\Bbb M}_*^n$ is
equicontinuous whenever $f$ lies in a ball $B_R,$ $f$ does not take
values from a fixed continuum $K\subset B_R,$ and a
quasiconformality coefficient $Q(x)$ has a finite mean oscillation
at every point.

\section{Введение} В недавно вышедшей работе \cite{IS} был получен
некоторый результат о равностепенной непрерывности одного класса
отображений между римановыми многообразиями. Речь идёт о локальном
поведении так называемых кольцевых $Q$-отображений (отображений с
неограниченной характеристикой квазиконформности $Q,$ искажающих
конформный модуль семейств кривых в $Q(x)$ раз, где $Q$ -- наперёд
заданная положительная функция). Основная цель настоящей заметки --
показать справедливость аналогичного утверждения не только для
<<конформного>> модуля $M,$ но и когда искажение семейств кривых
происходит относительно модуля порядка $p\in [n-1, n],$ где $n$ --
размерность риманова многообразия ${\Bbb M}^n,$ а $p$ --
фиксированное число из указанного отрезка. Стоит заметить, что
согласно Герингу в пространстве ${\Bbb R}^n$ такие отображения
квазиизометричны при ограниченных $Q,$ другими словами, при
некоторой постоянной $C>0$ и всех $x_0\in D$ справедлива оценка
%
$$\limsup\limits_{x\rightarrow
x_0}\frac{|f(x)-f(x_0)|}{|x-x_0|}\leqslant C\,,$$
%
см., напр., \cite[теорема~2]{Ge}. При $p=n$ свойство
квазиизометричности указанных отображений, к сожалению,
утрачивается, как показывает простой пример отображения с
ограниченным искажением $f(x)=x|x|^{\alpha-1},$ $x\in {\Bbb
R}^n\setminus\{0\},$ $0<\alpha<1,$ $f(0):=0.$

\medskip
Перейдём к определениям и формулировкам основных результатов.
Следующие понятия могут быть найдены, напр., в \cite{Lee} и
\cite{PSh}. Напомним, что {\it $n$-мерным топологическим
многообразием} ${\Bbb M}^n$ называется хаусдорфово топологическое
пространство со счётной базой, каждая точка которого имеет
окрестность, гомеоморфную некоторому открытому множеству в ${\Bbb
R}^n.$ {\it Картой} на многообразии ${\Bbb M}^n$ будем называть пару
$(U, \varphi),$ где $U$ --- открытое множество пространства ${\Bbb
M}^n$ и $\varphi$ --- соответствующий гомеоморфизм множества $U$ на
открытое множество в ${\Bbb R}^n.$ Если $p\in U$ и
$\varphi(p)=(x^1,\ldots,x^n)\in {\Bbb R}^n,$ то соответствующие
числа $x^1,\ldots,x^n$ называются {\it локальными координатами
точки} $p.$ {\it Гладким многообразием} называется само множество
${\Bbb M}^n$ вместе с соответствующим набором карт $(U_{\alpha},
\varphi_{\alpha}),$ так, что объединение всех $U_{\alpha}$ по
параметру $\alpha$ даёт всё ${\Bbb M}^n$ и, кроме того, отображение,
осуществляющее переход от одной системы локальных координат к
другой, принадлежит классу $C^{\infty}.$

\medskip
Напомним, что {\it римановой метрикой} на гладком многообразии
${\Bbb M}^n$ называется положительно определённое гладкое
симметричное тензорное поле типа $(0,2).$ В частности, компоненты
римановой метрики $g_{kl}$ в различных локальных координатах $(U,
x)$ и $(V, y)$ взаимосвязаны посредством тензорного закона
$$'g_{ij}(x)=g_{kl}(y(x))\frac{\partial y^k}{\partial x^i}
\frac{\partial y^l}{\partial x^j}.$$
{\it Римановым многообразием} будем называть гладкое многообразие
вместе с римановой метрикой на нём. Длину гладкой кривой
$\gamma=\gamma(t),$ $t\in [t_1, t_2],$ соединяющей точки
$\gamma(t_1)=M_1\in {\Bbb M}^n,$ $\gamma(t_2)=M_2\in {\Bbb M}^n$ и
$n$-мерный объём ({\it меру объёма $v$}) множества $A$ на римановом
многообразии определим согласно соотношениям
\begin{equation}\label{eq8}
l(\gamma):=\int\limits_{t_1}^{t_2}\sqrt{g_{ij}(x(t))\frac{dx^i}{dt}\frac{dx^j}{dt}}\,dt,\quad
v(A)=\int\limits_{A}\sqrt{\det g_{ij}}\,dx^1\ldots dx^n\,.
\end{equation}
Ввиду положительной определённости тензора $g=g_{ij}(x)$ имеем:
$\det g_{ij}>0.$ {\it Геодезическим расстоянием} между точками $p_1$
и $p_2\in {\Bbb M}^n$ будем называть наименьшую длину всех
кусочно-гладких кривых в ${\Bbb M}^n,$ соединяющих точки $p_1$ и
$p_2.$ Геодезическое расстояние между точками $p_1$ и $p_2$ будем
обозначать символом $d(p_1, p_2)$ (всюду далее $d$ обозначает
геодезическое расстояние, если не оговорено противное). Так как
риманово многообразие, вообще говоря, не предполагается связным,
расстояние между любыми точками многообразия, вообще говоря, может
быть не определено. Хорошо известно, что любая точка $p$ риманова
многообразия ${\Bbb M}^n$ имеет окрестность $U\ni p$ (называемую
далее {\it нормальной окрестностью точки $p$}) и соответствующее
координатное отображение $\varphi\colon U\rightarrow {\Bbb R}^n,$
так, что геодезические сферы с центром в точке $p$ и радиуса $r,$
лежащие в окрестности $U,$ переходят при отображении $\varphi$ в
евклидовы сферы того же радиуса, а пучок геодезических кривых,
исходящих из точки $p,$ переходит в пучок радиальных отрезков в
${\Bbb R}^n$ (см.~\cite[леммы~5.9 и 6.11]{Lee}, см.\ также
комментарии на стр.~77 здесь же). Локальные координаты
$\varphi(p)=(x^1,\ldots, x^n)$ в этом случае называются {\it
нормальными координатами} точки $p.$ Стоит отметить, что в случае
связного многообразия ${\Bbb M}^n$ открытые множества метрического
пространства $({\Bbb M}^n, d)$ порождают топологию исходного
топологического пространства ${\Bbb M}^n$
(см.~\cite[лемма~6.2]{Lee}). Заметим, что в нормальных координатах
всегда тензорная матрица $g_{ij}(x)$ в точке $p$ --- единичная (а в
силу непрерывности $g$ в точках, близких к $p,$ эта матрица сколь
угодно близка к единичной; см.~\cite[пункт~(c)
предложения~5.11]{Lee}).

\medskip
Семейство $\frak{F}$ отображений $f\colon X\rightarrow
{X}^{\,\prime}$ называется {\it равностепенно непрерывным в точке}
$x_0 \in X,$ если для любого $\varepsilon>0$ найдётся такое
$\delta>0$, что ${d}^{\,\prime}
\left(f(x),f(x_0)\right)<\varepsilon$ для всех таких $x,$ что
$d(x,x_0)<\delta$ и для всех $f\in \frak{F}.$ Говорят, что
$\frak{F}$ {\it равностепенно непрерывно}, если $\frak{F}$
равностепенно непрерывно в каждой  точке $x_0\in X.$ Согласно одной
из версий теоремы Арцела--Асколи (см., напр.,
\cite[пункт~20.4]{Va}), если $\left(X,\,d\right)$ --- сепарабельное
метрическое пространство, а $\left(X^{\,\prime},\,
d^{\,\prime}\right)$ --- компактное метрическое пространство, то
семейство $\frak{F}$ отображений $f\colon X\rightarrow
{X}^{\,\prime}$ нормально тогда  и только тогда, когда  $\frak{F}$
равностепенно непрерывно. Здесь и далее равностепенная непрерывность
семейства отображений $\{f\colon {\Bbb M}^n\rightarrow {\Bbb
M}_*^n\}$ понимается в смысле геодезических расстояний $d$ и
$d^{\,\prime}$ на римановых многообразиях ${\Bbb M}^n$ и ${\Bbb
M}_*^n,$ соответственно.

\medskip
Пусть $\left(X,\,d, \mu\right)$ --- произвольное метрическое
пространство, наделённое локально конечной борелевской мерой $\mu$ и
$$B(x_0, r)=\{x\in X: d(x, x_0)<r\}\,.$$
Следующее определение может быть найдено, напр., в
\cite[разд.~4]{RSa}. Будем говорить, что интегрируемая в $B(x_0, r)$
функция ${\varphi}\colon D\rightarrow{\Bbb R}$ имеет {\it конечное
среднее колебание} в точке $x_0\in D$, пишем $\varphi\in FMO(x_0),$
если
%
%
%
%
$$\limsup\limits_{\varepsilon\rightarrow 0}\frac{1}{\mu(B(
x_0,\,\varepsilon))}\int\limits_{B(x_0,\,\varepsilon)}
|{\varphi}(x)-\overline{{\varphi}}_{\varepsilon}|\,
d\mu(x)<\infty,$$
%
%
где
$\overline{{\varphi}}_{\varepsilon}=\frac{1} {\mu(B(
x_0,\,\varepsilon))}\int\limits_{B( x_0,\,\varepsilon)}
{\varphi}(x)\, d\mu(x).$

\medskip
Всюду далее (если не оговорено противное) ${\Bbb M}^n$ и ${\Bbb
M}_*^n$ -- римановы многообразия с геодезическими расстояниями $d$ и
$d_*,$ соответственно. {\it Кривой} $\gamma$ мы называем непрерывное
отображение отрезка $[a,b]$ (открытого интервала $(a,b),$ либо
полуоткрытого интервала вида $[a,b)$ или $(a,b]$) в ${\Bbb M}^n,$
$\gamma\colon [a,b]\rightarrow {\Bbb M}^n.$ Под семейством кривых
$\Gamma$ подразумевается некоторый фиксированный набор кривых
$\gamma,$ а, если $f\colon{\Bbb M}^n\rightarrow {\Bbb M}_*^n$ ---
произвольное отображение, то
$f(\Gamma)=\left\{f\circ\gamma|\gamma\in\Gamma\right\}.$ Длину
произвольной кривой $\gamma\colon [a, b]\rightarrow {\Bbb M}^n,$
лежащей на многообразии ${\Bbb M}^n,$ можно определить как точную
верхнюю грань сумм $\sum\limits_{i=1}^{n-1} d(\gamma(t_i),
\gamma(t_{i+1}))$ по всевозможным разбиениям $a\leqslant
t_1\leqslant\ldots\leqslant t_n\leqslant b.$ Следующие определения в
случае пространства ${\Bbb R}^n$ могут быть найдены, напр., в
\cite[разд.~1--6, гл.~I]{Va}, см.\ также \cite[гл.~I]{Fu}. Борелева
функция $\rho\colon {\Bbb M}^n\,\rightarrow [0,\infty]$ называется
{\it допустимой} для семейства $\Gamma$ кривых $\gamma$ в ${\Bbb
M}^n,$ если линейный интеграл по натуральному параметру $s$ каждой
(локально спрямляемой) кривой $\gamma\in \Gamma$ от функции $\rho$
удовлетворяет условию $\int\limits_{0}^{l(\gamma)}\rho
(\gamma(s))ds\geqslant 1.$ В этом случае мы пишем:
$\rho\in\mathrm{adm}\,\Gamma.$

\medskip
Зафиксируем $p\geqslant 1,$ тогда {\it $p$-мо\-ду\-лем} семейства
кривых $\Gamma $ называется величина
$$M_p(\Gamma)=\inf\limits_{\rho\in\mathrm{adm}\,\Gamma}
\int\limits_D \rho ^p (x)\,dv(x).$$
(Здесь и далее $v$ означает меру объёма, определённую в
(\ref{eq8})). При этом, если $\mathrm{adm}\,\Gamma=\varnothing,$ то
полагаем: $M_p(\Gamma)=\infty$ (см.~\cite[разд.~6 на с.~16]{Va} либо
\cite[с.~176]{Fu}). Свойства модуля в некоторой мере аналогичны
свойствам меры Лебега $m$ в ${\Bbb R}^n.$ Именно, модуль пустого
семейства кривых равен нулю, $M_p(\varnothing)=0,$ обладает
свойством
монотонности относительно семейств кривых, %
$ \Gamma_1\subset\Gamma_2\Rightarrow M_p(\Gamma_1)\leqslant
M_p(\Gamma_2), $
а также свойством полуаддитивности:
$M_p\left(\bigcup\limits_{i=1}^{\infty}\Gamma_i\right)\leqslant
\sum\limits_{i=1}^{\infty}M_p(\Gamma_i) $
(см.~\cite[теорема~6.2, гл.~I]{Va} в ${\Bbb R}^n$ либо
\cite[теорема~1]{Fu} в случае более общих пространств с мерами).
Говорят, что семейство кривых $\Gamma_1$ \index{минорирование}{\it
минорируется} семейством $\Gamma_2,$ пишем $\Gamma_1\,>\,\Gamma_2,$
если для каждой кривой $\gamma\,\in\,\Gamma_1$ существует подкривая,
которая принадлежит семейству $\Gamma_2.$
В этом случае,
\begin{equation}\label{eq32*A}
\Gamma_1
> \Gamma_2 \quad \Rightarrow \quad M_p(\Gamma_1)\leqslant M_p(\Gamma_2)
\end{equation} (см.~\cite[теорема~6.4, гл.~I]{Va} либо
\cite[свойство~(c)]{Fu} в случае более общих пространств с мерами).

\medskip{}
Следующее определение для случая ${\Bbb R}^n$ может быть найдено,
напр., в работе \cite{SS}. Пусть $D$ -- область в ${\Bbb M}^n,$
$x_0\in D,$ $Q\colon D\rightarrow [0,\infty]$ --- измеримая
относительно меры объёма $v$ функция, и число $r_0>0$ таково, что
замкнутый шар $\overline{B(x_0, r_0)}$ лежит в некоторой нормальной
окрестности $U$ точки $x_0.$ Пусть также $0<r_1<r_2<r_0,$
$$A=A(x_0, r_1,r_2)=\{x\in {\Bbb M}^n: r_1<d(x, x_0)<r_2\}\,,$$
$$S_i=S(x_0,r_i)=\{x\in {\Bbb M}^n: d(x, x_0)=r_i\}, \quad i=1,2\,,$$
--- геодезические сферы с центром в точке
$x_0$ и радиусов $r_1$ и $r_2,$ соответственно, а
$\Gamma\left(S_1,\,S_2,\,A\right)$ обозначает семейство всех кривых,
соединяющих $S_1$ и $S_2$ внутри области $A.$
%
Отображение $f\colon D\rightarrow {\Bbb M}_*^n$ условимся называть
{\it кольцевым $(p, Q)$-отображением в точке $x_0\,\in\,D,$} если
соотношение
%
$$M_p\left(f\left(\Gamma\left(S_1,\,S_2,\,A\right)\right)\right)\
\leqslant \int\limits_{A} Q(x)\cdot \eta^p(d(x, x_0))\ dv(x)$$
выполнено в кольце $A$ для произвольных $r_1,r_2,$ указанных выше, и
для каждой измеримой функции $\eta \colon  (r_1,r_2)\rightarrow
[0,\infty ]\,$ такой, что
$\int\limits_{r_1}^{r_2}\eta(r)dr\geqslant 1.$ Отображения типа
кольцевых $Q$-отображений были предложены к изучению О.~Мартио и
изучались им совместно с В.~Рязановым, У.~Сребро и Э.~Якубовым,
см.~\cite{MRSY}.

\medskip{}
Пусть $(X, d, \mu)$ --- метрическое пространство с метрикой $d,$
наделённое локально конечной борелевской мерой $\mu.$ Следуя
\cite[раздел 7.22]{He} будем говорить, что борелева функция
$\rho\colon  X\rightarrow [0, \infty]$ является {\it верхним
градиентом} функции $u\colon X\rightarrow {\Bbb R},$ если для всех
спрямляемых кривых $\gamma,$ соединяющих точки $x$ и $y\in X,$
выполняется неравенство $|u(x)-u(y)|\leqslant
\int\limits_{\gamma}\rho\,ds,$ где, как обычно,
$\int\limits_{\gamma}\rho\,ds$ обозначает линейный интеграл от
функции $\rho$ по кривой $\gamma.$ Будем также говорить, что в
указанном пространстве $X$ выполняется $(1; p)$-неравенство
Пуанкаре, если найдутся постоянные $C\geqslant 1$ и $\tau>0$ так,
что для каждого шара $B\subset X,$ произвольной ограниченной
непрерывной функции $u\colon B\rightarrow {\Bbb R}$ и любого её
верхнего градиента $\rho$ выполняется следующее неравенство:
$$\frac{1}{\mu(B)}\int\limits_{B}|u-u_B|d\mu(x)\leqslant C\cdot({\rm diam\,}B)\left(\frac{1}{\mu(\tau B)}
\int\limits_{\tau B}\rho^p d\mu(x)\right)^{1/p}\,,$$
где ${\rm diam}\,A$ --- диаметр множества $A\subset X,$ а
$u_B:=\frac{1}{\mu(B)}\int\limits_{B}u d\mu(x).$
Метрическое пространство $(X, d, \mu)$ назовём {\it
$\widetilde{Q}$-регулярным по Альфорсу} при некотором
$\widetilde{Q}\geqslant 1,$ если при каждом $x_0\in X,$ некоторой
постоянной $C\geqslant 1$ и произвольного $R<{\rm diam}\,X,$
$\frac{1}{C}R^{\widetilde{Q}}\leqslant \mu(B(x_0, R))\leqslant
CR^{\widetilde{Q}}.$
Заметим, что локально римановы многообразия являются $n$-регулярными
по Альфорсу (см.~\cite[лемма~5.1]{ARS}). Следует также заметить, что
если риманово многообразие $\widetilde{Q}$-регулярно по Альфорсу, то
$\widetilde{Q}=n$ (см.\ рассуждения на с.~61 в \cite{He} о
совпадении $\widetilde{Q}$ с хаусдорфовой размерностью пространства
$X,$ а также \cite[лемма~5.1]{ARS} о совпадении топологической и
хаусдорфовых размерностей областей риманового многообразия).
Справедлива следующая

\medskip
\begin{theorem}\label{theor4*!} {\,
Пусть $p\in [n-1, n]$ и $\delta>0,$ многообразие ${\Bbb M}_*^n$
связно{\em,} является $n$-ре\-гу\-ляр\-ным по Альфорсу{\em,} кроме
того{\em,} в ${\Bbb M}_*^n$ выполнено $(1;p)$-неравенство Пуанкаре.
Пусть $B_R\subset {\Bbb M}_*^n$ --- некоторый фиксированный шар
радиуса $R,$ $D$
--- область в ${\Bbb M}^n$ и $Q\colon D\rightarrow [0, \infty]$ ---
функция{\em,} измеримая относительно меры объёма $v.$ Обозначим
через $\frak{R}_{x_0, Q, B_R, \delta, p}(D)$ семейство открытых
дискретных кольцевых $(p, Q)$-отображений $f\colon D\rightarrow B_R$
в точке $x_0\in D,$ для которых существует континуум $K_f\subset
B_R$ такой, что $f(x)\not\in K_f$ при всех $x\in D$ и, кроме того,
${\rm diam}\, K_f\geqslant \delta.$ Тогда семейство отображений
$\frak{R}_{x_0, Q, B_R, \delta, p}(D)$ является равностепенно
непрерывным в точке $x_0\in D,$ если $Q\in FMO(x_0).$}
\end{theorem}

\medskip
Ввиду теоремы Арцела--Асколи имеем также следующее

\medskip
\begin{corollary}\label{cor1} {\, Предположим, что в условиях
теоремы~{\ref{theor4*!}} многообразие ${\Bbb M}_*^n$ является
компактным, тогда $\frak{R}_{x_0, Q, B_R, \delta, p}(D)$ образует
нормальное семейство отображений. }
\end{corollary}

\medskip
\section{Вспомогательные леммы} Всюду далее граница $\partial D$
области $D\subset {\Bbb M}^n$ и замыкание $\overline{D}$ области $D$
понимаются в смысле геодезического расстояния $d.$ Перед тем, как мы
приступим к изложению вспомогательных результатов и основной части
данного раздела, дадим ещё одно важное определение
(см.~\cite[раздел~3, гл.~II]{Ri}). Пусть $D$ --- область риманового
многообразия ${\Bbb M}^n,$ $n\geqslant 2,$ $f\colon D \rightarrow
{\Bbb M}_*^n$ --- отображение, $\beta\colon [a,\,b)\rightarrow {\Bbb
M}_*^n$ --- некоторая кривая и
$x\in\,f^{\,-1}\left(\beta(a)\right).$ Кривая $\alpha\colon
[a,\,c)\rightarrow D,$ $a<c\leqslant b,$ называется {\it
максимальным поднятием} кривой $\beta$ при отображении $f$ с началом
в точке $x,$ если $(1)\quad \alpha(a)=x;$ $(2)\quad
f\circ\alpha=\beta|_{[a,\,c)};$ $(3)$\quad если
$c<c^{\prime}\leqslant b,$ то не существует кривой
$\alpha^{\prime}\colon [a,\,c^{\prime})\rightarrow D,$ такой что
$\alpha=\alpha^{\prime}|_{[a,\,c)}$ и $f\circ
\alpha=\beta|_{[a,\,c^{\prime})}.$ Имеет место следующее

\medskip
\begin{proposition}\label{pr7}
{ Пусть $D$ --- область в ${\Bbb M}^n,$ $f\colon D\rightarrow {\Bbb
M}^n_*$ --- открытое дискретное отображение{\em,} $\beta\colon
[a,\,b)\rightarrow {\Bbb M}_*^n$ --- кривая и точка
$x\in\,f^{-1}\left(\beta(a)\right).$ Тогда кривая $\beta$ имеет
максимальное поднятие при отображении $f$ с началом в точке $x,$ см.
\cite[предложение~2.1]{IS}. }
\end{proposition}

\medskip{}
Пусть $A$ --- открытое подмножество многообразия ${\Bbb M}^n,$ а $C$
--- компактное подмножество $A.$  {\it Конденсатором} будем называть
пару множеств $E=\left(A,\,C\right).$ Пусть $p\geqslant 1,$ тогда
{\it $p$-ёмкостью} конденсатора $E$ будем называть следующую
величину: ${\rm cap}_p\,E=M_p(\Gamma_E),$ где $\Gamma_E$
--- семейство всех кривых вида $\gamma\colon [a,\,b)\rightarrow A,$
таких что $\gamma(a)\in C$ и $|\gamma|\cap\left(A\setminus
F\right)\ne\varnothing$ для произвольного компакта $F\subset A.$
Заметим, что в случае пространства ${\Bbb R}^n$ указанная величина
${\rm cap}_p\,E$ совпадает с классическим определением $p$-ёмкости
(см.~\cite[предложение~10.2 и замечание~10.8, гл.~II]{Ri}).

\medskip{}
Следующая лемма может быть полезной при исследовании свойства
равностепенной непрерывности открытых дискретных кольцевых $(p,
Q)$-отоб\-ра\-же\-ний в наиболее общей ситуации. Её доказательство
аналогично случаю ${\Bbb R}^n$ (см. \cite[лемма~4.1]{SS}), однако,
для полноты и строгости изложения мы приводим его полностью для
случая произвольных римановых многообразий.

\medskip
\begin{lemma}\label{lem4}{\,
Пусть $p\geqslant 1,$ $D$ --- область в ${\Bbb M}^n,$ $f\colon
D\rightarrow{\Bbb M}_*^n$ --- открытое дискретное кольцевое $(p,
Q)$-отоб\-ра\-же\-ние в точке $x_0\in D,$ $r_0>0$ таково{\em,} что
шар $B(x_0, r_0)$ лежит со своим замыканием в некоторой нормальной
окрестности $U$ точки $x_0.$ Предположим{\em,} что для некоторого
числа $0<\varepsilon_0<r_0,$ некоторого $\varepsilon_0^{\,\prime}\in
(0, \varepsilon_0)$ и семейства неотрицательных измеримых по Лебегу
функций $\{\psi_{\varepsilon}(t)\},$ $\psi_{\varepsilon}\colon
(\varepsilon, \varepsilon_0)\rightarrow [0, \infty],$
$\varepsilon\in\left(0, \varepsilon_0^{\,\prime}\right),$ выполнено
условие
 \begin{equation} \label{eq3.7B}
\int\limits_{\varepsilon<d(x, x_0)<\varepsilon_0}
Q(x)\cdot\psi_{\varepsilon}^p(d(x, x_0))\, dv(x)\leqslant
F(\varepsilon, \varepsilon_0)\qquad\forall\,\,\varepsilon\in(0,
\varepsilon_0^{\,\prime}),
 \end{equation}
где $F(\varepsilon, \varepsilon_0)$ --- некоторая функция и
\begin{equation}\label{eq3AB} 0<I(\varepsilon, \varepsilon_0):=
\int\limits_{\varepsilon}^{\varepsilon_0}\psi_{\varepsilon}(t)dt <
\infty\qquad\forall\,\,\varepsilon\in(0,
\varepsilon_0^{\,\prime}).\end{equation}
Тогда
\begin{equation}\label{eq3B}
{\rm cap}_p\,f(E)\leqslant F(\varepsilon,\varepsilon_0)/
I^{p}(\varepsilon, \varepsilon_0)\qquad
\forall\,\,\varepsilon\in\left(0,\,\varepsilon_0^{\,\prime}\right),
\end{equation}
где $E=\left(A,\,C\right)$ --- конденсатор{\em,} $A=B(x_0,
\varepsilon_0)$ и $C=\overline{B(x_0, \varepsilon)}.$
}
\end{lemma}

\medskip
\begin{proof}
Рассмотрим конденсатор $E=(A,\,C),$ где $A$ и $C$ таковы, как
указано в условии леммы. Если ${\rm cap}_p\,f(E)=0,$ доказывать
нечего. Пусть ${\rm cap}_p\,f(E)\ne 0.$

\medskip
Пусть $\Gamma_E$ --- семейство всех кривых, соответствующее
определению $p$-ёмкости. Для конденсатора $f(E)$ рассмотрим также
семейство кривых $\Gamma_{f(E)}.$ Заметим также, что каждая кривая
$\gamma\in\Gamma_{f(E)}$ имеет максимальное поднятие при отображении
$f,$ лежащее в $A$ с началом в $C$ (см.\ предложение~\ref{pr7}).
Пусть $\Gamma^{\,*}$ --- семейство всех максимальных под\-ня\-тий
кривых $\Gamma_{f(E)}$ при отображении $f$ с началом в $C.$ Заметим,
что $\Gamma^{*}\subset \Gamma_E$ (см. подробности в доказательстве
\cite[лемма~2.1]{IS}).

\medskip
Заметим, что $\Gamma_{f(E)}>f(\Gamma^{*}),$ и, следовательно, ввиду
свойства~\eqref{eq32*A}
\begin{equation}\label{eq7AB}
M_p\left(\Gamma_{f(E)}\right)\leqslant
M_p\left(f(\Gamma^{*})\right).
\end{equation}
Рассмотрим
$S_{\,\varepsilon}=S(x_0,\,\varepsilon) = \{ x\,\in\,{\Bbb M}^n :
d(x, x_0)=\,\varepsilon\},$ $
S_{\,\varepsilon_{0}}=S(x_0,\,\varepsilon_0) = \{ x\,\in\,{\Bbb M}^n
: d(x, x_0)=\,\varepsilon_0\},$
где $\varepsilon_0$ --- из условия леммы и
$\varepsilon\in\left(0,\,\varepsilon_0^{\,\prime}\right).$ Заметим,
что, поскольку $\Gamma^{*}\,\subset\,\Gamma_E,$ то при всех
достаточно малых $\delta>0$ будем иметь, что
$\Gamma\left(S_{\varepsilon}, S_{\varepsilon_0-\delta}, A(x_0,
\varepsilon, \varepsilon_0-\delta)\right)<\Gamma^{*}$ и,
следовательно, $f(\Gamma\left(S_{\varepsilon},
S_{\varepsilon_0-\delta}, A(x_0, \varepsilon,
\varepsilon_0-\delta)\right))<f(\Gamma^{\,*}).$ Здесь мы положили
$$S_{\varepsilon_0-\delta}:=\{x\in {\Bbb M}^n:
d(x, x_0)<\varepsilon_0-\delta\}\,.$$ Значит,
\begin{equation}\label{eq3.3.1}
M_p\left(f(\Gamma^{*})\right)\leqslant
 M_p\left(f\left(\Gamma\left(S_{\varepsilon}, S_{\varepsilon_0-\delta},
 A(x_0, \varepsilon, \varepsilon_0-\delta)\right)\right)\right).
\end{equation}
Из соотношений~\eqref{eq7AB} и~\eqref{eq3.3.1} следует, что
$ M_p\left(\Gamma_{f(E)}\right)\,\leqslant
 M_p\left(f\left(\Gamma\left(S_{\varepsilon}, S_{\varepsilon_0-\delta},
 A(x_0, \varepsilon, \varepsilon_0-\delta)\right)\right)\right)
$
и, таким образом, по определению
\begin{equation}\label{eq5aa}
{\rm cap}_p\,\,f(E) \leqslant
 M_p\left(f\left(\Gamma\left(S_{\varepsilon}, S_{\varepsilon_0-\delta},
 A(x_0, \varepsilon, \varepsilon_0-\delta)\right)\right)\right).
\end{equation}
Пусть $\eta(t)$ произвольная неотрицательная измеримая функция,
удовлетворяющая условию
$\int\limits_{\varepsilon}^{\varepsilon_0}\eta(t)dt=1.$ Рассмотрим
семейство измеримых функций
$\eta_{\delta}(t)=\frac{\eta(t)}{\int\limits_{\varepsilon}^{\varepsilon_0-\delta}\eta(t)dt}.$
(Так как $\int\limits_{\varepsilon}^{\varepsilon_0}\eta(t)dt=1,$ то
$\delta>0$ можно выбрать так, что
$\int\limits_{\varepsilon}^{\varepsilon_0-\delta}\eta(t)dt>0$).
Поскольку
$\int\limits_{\varepsilon}^{\varepsilon_0-\delta}\eta_{\delta}(t)dt=1,$
то по определению кольцевого $(p, Q)$-отображения в точке $x_0$ мы
получим
 \begin{multline*}
M_p\left(f\left(\Gamma\left(S_{\varepsilon}\,,
S_{\varepsilon_0-\delta}, A(x_0, \varepsilon, \varepsilon_0-\delta)\right)\right)\right)\,\leqslant\\
\leqslant\frac{1}{\left(\int\limits_{\varepsilon}^{\varepsilon_0-\delta}\eta(t)dt\right)^p}
\int\limits_{\varepsilon<d(x,
x_0)<\varepsilon_0}Q(x)\cdot\eta^p(d(x, x_0))\, \ dv(x).
 \end{multline*}
Переходя здесь к пределу при $\delta\rightarrow 0$ и учитывая
соотношение \eqref{eq5aa}, получаем что
\begin{equation}\label{eq1F}
{\rm cap}_p\,\,f(E) \leqslant \int\limits_{\varepsilon<d(x,
x_0)<\varepsilon_0}Q(x)\cdot\eta^p(d(x, x_0))\, \
dv(x)\end{equation}
для произвольной неотрицательной измеримой функции $\eta(t)$ такой
что $\int\limits_{\varepsilon}^{\varepsilon_0}\eta(t)dt=1.$
Рассмотрим семейство измеримых функций
$\eta_{\varepsilon}(t)=\psi_{\varepsilon}(t)/I(\varepsilon,
\varepsilon_0 ),$ $t\in(\varepsilon,\, \varepsilon_0).$ Заметим, что
для $\varepsilon\in (0, \varepsilon_0^{\,\prime})$
выполнено равенство
$\int\limits_{\varepsilon}^{\varepsilon_0}\eta_{\varepsilon}(t)\,dt=1.$
Тогда из (\ref{eq1F}) мы получим, что для любого $\varepsilon\in (0,
\varepsilon_0^{\,\prime})$
 \begin{equation}\label{eq8A}
{\rm cap}_p\,\,f(E)\leqslant\frac{1}{I^p(\varepsilon,
\varepsilon_0)} \int\limits_{\varepsilon<d(x,
x_0)<\varepsilon_0}Q(x)\cdot\psi_{\varepsilon}^p(d(x, x_0))\,dv(x).
 \end{equation}
Из соотношений~\eqref{eq3.7B} и \eqref{eq8A} следует
соотношение~\eqref{eq3B}.~\end{proof}

\medskip
Сравнительно недавно Адамовичем и Шанмугалингам получено следующее
фундаментальное утверждение, которым нам необходимо воспользоваться
(см.~\cite[предложение~4.7]{AS}).

\medskip
\begin{proposition}\label{pr2}
{ Пусть $X$ --- $\widetilde{Q}$-регулярное по Альфорсу метрическое
пространство с мерой{\em,} в котором выполняется $(1;p)$-неравенство
Пуанкаре так{\em,} что $\widetilde{Q}-1\leqslant n\leqslant
\widetilde{Q}.$ Тогда для произвольных континуумов $E$ и $F,$
содержащихся в шаре $B(x_0, R),$ и некоторой постоянной $C>0$
выполняется неравенство
$M_p(\Gamma(E, F, X))\geqslant \frac{1}{C}\cdot\frac{\min\{{\rm
diam}\,E, {\rm diam}\,F\}}{R^{1+p-\widetilde{Q}}}.$ }
\end{proposition}

\medskip
Следующая лемма является ключом к доказательству основных
результатов работы.

\medskip
\begin{lemma}\label{lem1}{\,Пусть $p\in [n-1, n]$ и $\delta>0,$ многообразие ${\Bbb M}_*^n$
связно{\em,} является $n$-ре\-гу\-ляр\-ным по Альфорсу{\em,} кроме
того{\em,} в ${\Bbb M}_*^n$ выполнено $(1;p)$-неравенство Пуанкаре.
Пусть $B_R\subset {\Bbb M}_*^n$ --- некоторый фиксированный шар
радиуса $R,$ $D$
--- область в ${\Bbb M}^n$ и $Q\colon D\rightarrow [0, \infty]$ ---
функция{\em,} измеримая относительно меры объёма $v.$ Обозначим
через $\frak{R}_{x_0, Q, B_R, \delta, p}(D)$ семейство открытых
дискретных кольцевых $(p, Q)$-отображений $f\colon D\rightarrow B_R$
в точке $x_0\in D,$ для которых существует континуум $K_f\subset
B_R$ такой, что $f(x)\not\in K_f$ при всех $x\in D$ и, кроме того,
${\rm diam}\, K_f\geqslant \delta.$

Пусть $\varepsilon_0<r_0,$ где шар $B(x_0, r_0)$ лежит вместе со
своим замыканием в некоторой нормальной окрестности $U$ точки $x_0.$

Предположим также{\em,} что для некоторого числа
$\varepsilon_0^{\,\prime}\in (0, \varepsilon_0)$ и семейства
неотрицательных измеримых по Лебегу функций
$\{\psi_{\varepsilon}(t)\},$ $\psi_{\varepsilon}\colon (\varepsilon,
\varepsilon_0)\rightarrow [0, \infty],$ $\varepsilon\in\left(0,
\varepsilon_0^{\,\prime}\right),$ выполнено
условие~\eqref{eq3.7B}{\em,} где некоторая заданная функция
$F(\varepsilon, \varepsilon_0)$ удовлетворяет условию
$F(\varepsilon, \varepsilon_0)=o(I^n(\varepsilon, \varepsilon_0)),$
а $I(\varepsilon, \varepsilon_0)$ определяется
соотношением~\eqref{eq3AB}.

Тогда семейство отображений $\frak{R}_{x_0, Q, B_R, \delta, p}(D)$
является равностепенно непрерывным в точке $x_0\in D.$ }
\end{lemma}

\begin{proof}
Пусть $f\in\frak{R}_{x_0, Q, B_R, \delta, p}(D).$ Полагаем
$A:=B(x_0, r_0)\subset D.$ Заметим, что при указанных условиях
$\overline{A}$ является компактным подмножеством $D.$ Тогда при
каждом $0<\varepsilon<r_0$ множество $C:=\overline{B(x_0,
\varepsilon)}$ является компактным подмножеством $B(x_0, r_0).$
Таким образом, $E=(A, C)$ --- конденсатор в ${\Bbb M}^n.$

Рассмотрим семейство кривых $\Gamma_{f(E)}$ для конденсатора $f(E)$
в терминах определения $p$-ёмкости. Заметим, что подсемейство
неспрямляемых кривых семейства $\Gamma_{f(E)}$ имеет нулевой
$p$-модуль.

Действительно, $f(\overline{A})$ является компактным множеством
${\Bbb M}^n_*$ как непрерывный образ компакта $\overline{A}$ и,
значит, множество $f(A)$ имеет конечный объём $v_*(f(A))<\infty.$
Пусть $\Gamma^{\infty}_{f(E)}$ состоит из всех неспрямляемых кривых
семейства $\Gamma_{f(E)},$ тогда функция $\rho(x)=\varepsilon$
принадлежит ${\rm adm}\,\Gamma^{\infty}_{f(E)}.$ Ввиду этого, по
определению $p$-модуля, $M_p(\Gamma^{\infty}_{f(E)})\leqslant
\varepsilon^p\cdot v_*(f(A)).$ Устремляя здесь $\varepsilon$ к нулю,
получаем $M_p(\Gamma^{\infty}_{f(E)})=0,$ что и требовалось
установить.

Заметим также, что оставшееся подсемейство, состоящее из всех
спрямляемых кривых семейства $\Gamma_{f(E)},$ состоит из кривых
$\beta\colon [a, b)\rightarrow f(D),$ имеющих предел при
$t\rightarrow b.$ Заметим, что указанный предел принадлежит
множеству $\partial f(A).$ Из сказанного следует, что
\begin{equation}\label{eq1}
M_p(\Gamma_{f(E)})=M_p(\Gamma(f(C), \partial f(A), f(A))).
\end{equation}
Пусть $K_f$ -- тот континуум, который выпускает отображение $f.$
Заметим, что
$$\Gamma(K_f, f(C), {\Bbb M}_*^n)>\Gamma(f(C),
\partial f(A), f(A))$$
(см. \cite[теорема 1, $\S\,46,$ п.~I]{Ku}), так что ввиду
(\ref{eq32*A})
\begin{equation}\label{eq9}
M_p(\Gamma(f(C), \partial f(A), f(A)))\geqslant M_p(\Gamma(K_f,
f(C), {\Bbb M}_*^n))\,.
\end{equation}
Ввиду предложения \ref{pr2} получим:
\begin{equation}\label{eq2}
M_p(\Gamma(K_f, f(C), {\Bbb M}_*^n))\geqslant
\frac{1}{C_1}\cdot\frac{\min\{{\rm diam}\,f(C), {\rm
diam}\,K_f\}}{R^{1+p-n}}\,.
 \end{equation}
По лемме~\ref{lem4} $M_p(\Gamma_{f(E)})\rightarrow 0$ при
$\varepsilon\rightarrow 0,$ так что ввиду~\eqref{eq1} и~\eqref{eq9}
получаем, что
$$\min\{{\rm diam}\,f(C), {\rm diam}\,K_f\}={\rm diam}\,f(C)$$ при
$\varepsilon\rightarrow 0.$ Из соотношений~\eqref{eq3B}
и~\eqref{eq2} также вытекает, что для любого $\sigma>0$ найдётся
$\varepsilon_0=\varepsilon_0(\sigma)$ так, что при $\varepsilon\in
(0, \varepsilon_0)$
$${\rm diam}\,f(C)\leqslant \sigma,$$
что и означает равностепенную непрерывность семейства
$\frak{R}_{x_0, Q, B_R, \delta, p}(D)$ в точке $x_0.$~\end{proof}

\medskip
\section{Доказательство основных результатов} {\it Доказательство
теоремы~{\em\ref{theor4*!}}} вытекает из леммы~{\em\ref{lem1}}.
Выберем в этой лемме $0\,<\,\psi(t)\,=\,\frac
{1}{\left(t\,\log{\frac1t}\right)^{n/p}}.$ На основании
\cite[предложение~3]{Af$_1$} для указанной функции
будем иметь, что %
$$\int\limits_{\varepsilon<d(x, x_0)<\varepsilon_0}
Q(x)\cdot\psi^p(d(x, x_0))
 \ dv(x)\,= \int\limits_{\varepsilon<d(x, x_0)< {\varepsilon_0}}\frac{Q(x)\,
dv(x)} {\left(d(x, x_0) \log \frac{1}{d(x, x_0)}\right)^n} = $$
$$=O
\left(\log\log \frac{1}{\varepsilon}\right)$$
%
при  $\varepsilon \rightarrow 0.$
Заметим также, что при указанных выше $\varepsilon$ выполнено
$\psi(t)\geqslant \frac {1}{t\,\log{\frac1t}},$ поэтому
$I(\varepsilon,
\varepsilon_0)\,:=\,\int\limits_{\varepsilon}^{\varepsilon_0}\psi(t)\,dt\,\geqslant
\log{\frac{\log{\frac{1}
{\varepsilon}}}{\log{\frac{1}{\varepsilon_0}}}}.$ Тогда соотношение
(\ref{eq3.7B}) выполнено при $p=1.$ Таким образом, все условия леммы
\ref{lem1} выполнены и, значит, необходимое заключение вытекает из
этой леммы.~$\Box$

\medskip{} {\it Элементом площади} гладкой поверхности $H$ на
римановом многообразии ${\Bbb M}^n$ будем называть выражение вида
$d\mathcal{A}=\sqrt{{\rm det}\,g_{\alpha\beta}^*}\,du^1\ldots
du^{n-1},$
где $g_{\alpha\beta}^*$ --- риманова метрика на $H,$ порождённая
исходной римановой метрикой $g_{ij}$ согласно соотношению
\begin{equation}\label{eq5}
g_{\alpha\beta}^*(u)=g_{ij}(x(u))\frac{\partial x^i}{\partial
u^{\alpha}} \frac{\partial x^j}{\partial u^{\beta}}.
\end{equation}
Здесь индексы $\alpha$ и $\beta$ меняются от $1$ до $n-1,$ а $x(u)$
обозначает параметризацию поверхности $H$ такую, что $\nabla_u x\ne
0.$ Справедливо следующее утверждение, обобщающее теорему
\ref{theor4*!}.

\medskip
\begin{theorem}\label{th1} {\,
Пусть $p\in [n-1, n]$ и $\delta>0,$ многообразие ${\Bbb M}_*^n$
связно{\em,} является $n$-ре\-гу\-ляр\-ным по Альфорсу{\em,} кроме
того{\em,} в ${\Bbb M}_*^n$ выполнено $(1;p)$-неравенство Пуанкаре.
Пусть $B_R\subset {\Bbb M}_*^n$ --- некоторый фиксированный шар
радиуса $R,$ $D$
--- область в ${\Bbb M}^n$ и $Q\colon D\rightarrow [1, \infty]$ ---
функция{\em,} измеримая относительно меры объёма $v.$ Обозначим
через $\frak{R}_{x_0, Q, B_R, \delta, p}(D)$ семейство открытых
дискретных кольцевых $(p, Q)$-отображений $f\colon D\rightarrow B_R$
в точке $x_0\in D,$ для которых существует континуум $K_f\subset
B_R$ такой, что $f(x)\not\in K_f$ при всех $x\in D$ и, кроме того,
${\rm diam}\, K_f\geqslant \delta.$ Тогда семейство отображений
$\frak{R}_{x_0, Q, B_R, \delta, p}(D)$ является равностепенно
непрерывным в точке $x_0\in D,$ если при некотором $\delta(x_0)>0$
выполняется равенство
\begin{equation}\label{eq3}
\int\limits_{0}^{\delta(x_0)}\
\frac{dr}{r^{\frac{n-1}{p-1}}q_{x_0}^{\frac{1}{p-1}}(r)}=\infty\,,
\end{equation}
где $q_{x_0}(r):=\frac{1}{r^{n-1}}\int\limits_{S(x_0,
r)}Q(x)\,d\mathcal{A}.$ }
\end{theorem}

\begin{proof} Достаточно показать, что условие~\eqref{eq3} влечёт
выполнение соотношения~\eqref{eq3.7B} леммы~\ref{lem1}. Можно
считать, что $B(x_0, \delta(x_0))$ лежит в нормальной окрестности
точки $x_0.$ Рассмотрим функцию
\begin{equation}\label{eq1E} \psi(t)= \left \{\begin{array}{rr}
1/[t^{\frac{n-1}{p-1}}q_{x_0}^{\frac{1}{p-1}}(t)]\ , & \ t\in
(r_1,r_2)\ ,
\\ 0\ ,  &  \ t\notin (r_1,r_2)\ .
\end{array} \right.
\end{equation}
Заметим теперь, что требование вида~\eqref{eq3AB} выполняется при
$\varepsilon_0=\delta(x_0)$ и всех достаточно малых $\varepsilon.$
Далее установим неравенство
\begin{equation}\label{eq4}
\int\limits_{\varepsilon<d(x, x_0)<\delta(x_0)} Q(x)\psi^p(d(x,
x_0))\,dv(x)\leqslant C\cdot\int\limits_{\varepsilon}^{\delta(x_0)}
\frac{dt}{t^{\frac{n-1}{p-1}}q_{x_0}^{\frac{1}{p-1}}(t)}
\end{equation}
при некоторой постоянной $C>0.$ Для этого покажем, что к левой части
соотношения \eqref{eq4} применим аналог теоремы Фубини. Рассмотрим в
окрестности точки $x_0\in S(z_0, r)\subset {\Bbb R}^n$ локальную
систему координат $z^1,\ldots, z^n,$ $n-1$ базисных векторов которой
взаимно ортогональны и лежат в плоскости, касательной к сфере в
точке $x_0,$ а последний базисный вектор перпендикулярен этой
плоскости. Пусть $r, \theta^1,\ldots, \theta^{n-1}$ сферические
координаты точки $x=x(\theta)$ в ${\Bbb R}^n.$ Заметим, что $n-1$
приращений переменных $z^1,\ldots, z^{n-1}$ вдоль сферы при
фиксированном $r$ равны $dz^1=rd\theta^1,\dots,
dz^{n-1}=rd\theta^{n-1},$ а приращение переменной $z^n$ по $r$ равно
$dz^n=dr.$ В таком случае,
$$dv(x)=\sqrt{{\rm det\,}g_{ij}(x)}r^{n-1}\,dr d\theta^1\dots
d\theta^{n-1}.$$
Рассмотрим параметризацию сферы $S(0, r)$ $x=x(\theta),$
$\theta=(\theta^1,\ldots,\theta^{n-1}),$ $\theta_i\in (-\pi, \pi].$
Заметим, что $\frac{\partial x^{\alpha}}{\partial \theta^{\beta}}=1$
при $\alpha=\beta$ и $\frac{\partial x^{\alpha}}{\partial
\theta^{\beta}}=0$ при $\alpha\ne \beta,$ $\alpha,\beta=1,\ldots,
n-1.$ Тогда в обозначениях соотношения~\eqref{eq5} имеем:
$g_{\alpha\beta}^*(\theta)=g_{\alpha\beta}(x(\theta))r^2,$
$$d\mathcal{A}=\sqrt{\det\, g_{\alpha\beta}(x(\theta))}r^{n-1}d{\theta}^1\ldots d{\theta}^{n-1}.$$
Заметим, что
$$\frac{1}{r^{\frac{n-1}{p-1}}q_{x_0}^{\frac{1}{p-1}}(r)}=$$
\begin{equation}\label{eq6}
=\int\limits_{S(x_0,
r)}Q(x)\psi^p(d(x,x_0))\,d\mathcal{A}=\psi^p(r)r^{n-1}\cdot\int\limits_{\Pi}\sqrt{{\rm
det\,}g_{\alpha\beta}(x(\theta))}Q(x(\theta))\,d\theta^{1}\dots
d\theta^{n-1},
 \end{equation}
где $\Pi=(-\pi, \pi]^{n-1}$ --- прямоугольная область изменения
параметров $\theta^1,\ldots,\theta^{n-1}.$ Напомним, что в
нормальной системе координат геодезические сферы переходят в обычные
сферы того же радиуса с центром в нуле, а пучок геодезических,
исходящих из точки многообразия, переходит в пучок радиальных
отрезков в ${\Bbb R}^n$ (см.~\cite[леммы~5.9 и 6.11]{Lee}), так что
кольцу $\{x\in {\Bbb M}^n: \varepsilon<d(x, x_0)<\delta(x_0)\}$
соответствует та часть ${\Bbb R}^n,$ в которой $r\in (\varepsilon,
\delta(x_0)).$ Согласно сказанному выше, применяя классическую
теорему Фубини (см., напр.,~\cite[разд.~8.1, гл.~III]{Sa}),
 \begin{multline}\label{eq7}
\int\limits_{\varepsilon<d(x, x_0)<\delta(x_0)} Q(x)\psi^p(d(x,
x_0))\,dv(x)=\\
=\int\limits_{\varepsilon}^{\delta(x_0)}\int\limits_{\Pi} \sqrt{{\rm
det\,}g_{ij}(x)}Q(x)\psi^p(r)r^{n-1}\,d\theta^1\dots
d\theta^{n-1}dr.
 \end{multline}
Поскольку в нормальных координатах тензорная матрица $g_{ij}$  сколь
угодно близка к единичной в окрестности данной точки, то
$$C_2\det\,g_{\alpha\beta}(x)\leqslant\det\,g_{ij}(x)\leqslant
C_1\det\,g_{\alpha\beta}(x)\,.$$ Учитывая сказанное и сравнивая
\eqref{eq6} и \eqref{eq7}, приходим к соотношению~\eqref{eq4}. Но
тогда также
$$\int\limits_{\varepsilon<d(x, x_0)<\delta(x_0)}
Q(x)\psi^p(d(x, x_0))dv(x)=o(I^p(\varepsilon, \delta(x_0)))$$
ввиду соотношения~\eqref{eq3}. Утверждение теоремы следует теперь из
леммы~\ref{lem1}.~\end{proof}

\medskip
Ввиду теоремы Арцела--Асколи имеем также следующее

\medskip
\begin{corollary}\label{cor2} {\, Предположим{\em,} что в условиях
теоремы~{\em\ref{theor4*!}} многообразие ${\Bbb M}_*^n$ является
компактным. Тогда класс $\frak{R}_{x_0, Q, B_R, \delta, p}(D)$
образует нормальное семейство отображений. }
 \end{corollary}

\medskip
При $p\in [n-1, n)$ полученные выше результаты можно несколько
усилить.

\medskip
\medskip
\begin{theorem}\label{th2} {\,
Пусть $p\in [n-1, n)$ и $\delta>0,$ многообразие ${\Bbb M}_*^n$
связно{\em,} является $n$-ре\-гу\-ляр\-ным по Альфорсу{\em,} кроме
того{\em,} в ${\Bbb M}_*^n$ выполнено $(1;p)$-неравенство Пуанкаре.
Пусть $B_R\subset {\Bbb M}_*^n$ --- некоторый фиксированный шар
радиуса $R,$ $D$
--- область в ${\Bbb M}^n$ и $Q\colon D\rightarrow [0, \infty]$ ---
функция{\em,} измеримая относительно меры объёма $v.$ Обозначим
через $\frak{R}_{x_0, Q, B_R, \delta, p}(D)$ семейство открытых
дискретных кольцевых $(p, Q)$-отображений $f\colon D\rightarrow B_R$
в точке $x_0\in D,$ для которых существует континуум $K_f\subset
B_R$ такой, что $f(x)\not\in K_f$ при всех $x\in D$ и, кроме того,
${\rm diam}\, K_f\geqslant \delta.$ Тогда семейство отображений
$\frak{R}_{x_0, Q, B_R, \delta, p}(D)$ является равностепенно
непрерывным в точке $x_0\in D,$ если $Q\in L_{loc}^s({\Bbb R}^n)$
при некотором $s\geqslant\frac{n}{n-p}.$}
\end{theorem}

\medskip
\begin{proof}
Зафиксируем произвольным образом $0<\varepsilon_0<\infty,$ так,
чтобы шар $G:=B(x_0, \varepsilon_0)$ лежал вместе со своим
замыканием в некоторой нормальной окрестности точки $x_0,$ и положим
$\psi(t):=1/t.$ Заметим, что указанная функция $\psi$ удовлетворяет
неравенствам $0< I(\varepsilon,
\varepsilon_0):=\int\limits_{\varepsilon}^{\varepsilon_0}\psi(t)dt <
\infty.$ Покажем также, что в этом случае выполнено соотношение
\begin{equation} \label{eq4!}
\int\limits_{A(x_0,\varepsilon,
\varepsilon_0)}Q(x)\cdot\psi^{p}(d(x, x_0)) \
dv(x)\,=\,o\left(I^{p}(\varepsilon, \varepsilon_0)\right)\,.
\end{equation}
Применяя неравенство Гёльдера, будем иметь
$$\int\limits_{\varepsilon<d(x,
x_0)<\varepsilon_0}\frac{Q(x)}{d^{p}(x, x_0)} \ dv(x)\leqslant$$
\begin{equation}\label{eq13D}
\leqslant \left(\int\limits_{\varepsilon<d(x,
x_0)<\varepsilon_0}\frac{1}{d^{pq}(x, x_0)} \ dv(x)
\right)^{\frac{1}{q}}\,\left(\int\limits_{G} Q^{q^{\prime}}(x)\
dv(x)\right)^{\frac{1}{q^{\prime}}}\,,
\end{equation}
где  $1/q+1/q^{\prime}=1$. Заметим, что первый интеграл в правой
части неравенства (\ref{eq13D}) с точностью до некоторой постоянной
может быть подсчитан непосредственно. Действительно, пусть для
начала $q^{\prime}=\frac{n}{n-p}$ (и, следовательно,
$q=\frac{n}{p}.$) Ввиду аналога теоремы Фубини (см. ход
доказательства теоремы \ref{th1}) будем иметь:
$$
\int\limits_{\varepsilon<d(x, x_0)<\varepsilon_0}\frac{1}{d^{p q}(x,
x_0)} \ dv(x)\leqslant C\cdot
\int\limits_{\varepsilon}^{\varepsilon_0}
\frac{dt}{t}=C\cdot\log\frac{\varepsilon_0}{\varepsilon}\,.
$$
В обозначениях леммы \ref{lem4} мы будем иметь, что при
$\varepsilon\rightarrow 0$
$$
\frac{1}{I^{p}(\varepsilon,
\varepsilon_0)}\int\limits_{\varepsilon<d(x,
x_0)<\varepsilon_0}\frac{Q(x)}{d^{p}(x, x_0)} \ dv(x)\leqslant
C^{\frac{p}{n}}\Vert
Q\Vert_{L^{\frac{n}{n-p}}(G)}\left(\log\frac{\varepsilon_0}{\varepsilon}\right)
^{-p+\frac{p}{n}}\,\rightarrow 0\,,
$$
что влечёт выполнение соотношения (\ref{eq4!}).

\medskip
Пусть теперь $q^{\prime}>\frac{n}{n-p}$ (т.е.,
$q=\frac{q^{\prime}}{q^{\prime}-1}$). В этом случае
$$
\int\limits_{\varepsilon<d(x, x_0)<\varepsilon_0}\frac{1}{d^{pq}(x,
x_0)} \ dv(x)\leqslant C\int\limits_{\varepsilon}^{\varepsilon_0}
t^{n-\frac{p q^{\prime}}{q^{\prime}-1}-1}dt\leqslant
C\int\limits_{0}^{\varepsilon_0} t^{n-\frac{p
q^{\prime}}{q^{\prime}-1}-1}dt=$$
$$
=\frac{C}{n-\frac{p q^{\prime}}{q^{\prime}-1}}\varepsilon^{n-\frac{p
q^{\prime}}{q^{\prime}-1}}_0,
$$
и, значит,
$$
\frac{1}{I^{p}(\varepsilon, \varepsilon_0)}
\int\limits_{\varepsilon<d(x,
x_0)<\varepsilon_0}\frac{Q(x)}{d^{p}(x, x_0)} \ dv(x)\leqslant \Vert
Q\Vert_{L^{q^{\prime}}(G)}\left(\frac{C}{n-\frac{p
q^{\prime}}{q^{\prime}-1}} \varepsilon^{n-\frac{p
q^{\prime}}{q^{\prime}-1}}_0\right)^{\frac{1}{q}}\left(\log\frac{\varepsilon_0}{\varepsilon}\right)
^{-p}\,,
$$
что также влечёт выполнение соотношения (\ref{eq4!}).
Заметим, что оба соотношения (\ref{eq3.7B})--(\ref{eq3AB})
выполняются, так что желанное заключение вытекает из леммы
\ref{lem1}.~\end{proof}

{\footnotesize

? 

 }

\end{document}